\pdfoutput=1
\documentclass[11pt,a4paper]{article}

\usepackage[a4paper,margin=1in]{geometry}
\usepackage[T1]{fontenc}
\usepackage{lmodern}
\usepackage{microtype}
\usepackage{amsmath,amssymb,amsthm}
\usepackage{mathtools}
\usepackage{graphicx}
\usepackage{booktabs}
\usepackage{array}
\usepackage{enumitem}
\usepackage{xcolor}
\usepackage[font=small,labelfont=bf,labelsep=period]{caption}
\usepackage{authblk}
\usepackage{natbib}
\bibpunct[, ]{(}{)}{;}{a}{,}{,}

\usepackage[colorlinks=true,linkcolor=blue!45!black,citecolor=blue!45!black,urlcolor=blue!45!black,hyperfootnotes=false]{hyperref}
\hypersetup{
  pdftitle={Sparse optimal control for a class of infinite-dimensional bilinear systems: optimality conditions and uniqueness},
  pdfauthor={Redouane El Mezegueldy, Zakarya Dardour},
  pdfkeywords={Bilinear control; sparse optimal control; uniqueness of optimal controls; second-order optimality conditions}
}

\numberwithin{equation}{section}

\theoremstyle{plain}
\newtheorem{theorem}{Theorem}[section]
\newtheorem{lemma}[theorem]{Lemma}
\newtheorem{proposition}[theorem]{Proposition}
\newtheorem{corollary}[theorem]{Corollary}
\theoremstyle{definition}
\newtheorem{example}[theorem]{Example}
\theoremstyle{remark}
\newtheorem{remark}[theorem]{Remark}

\newcommand{\R}{\mathbb{R}}
\newcommand{\Uad}{U_{\mathrm{ad}}}
\newcommand{\dd}{\,\mathrm{d}}
\newcommand{\norm}[1]{\lVert #1\rVert}
\newcommand{\abs}[1]{\lvert #1\rvert}
\newcommand{\ip}[2]{\langle #1,\, #2\rangle}
\DeclareMathOperator{\supp}{supp}
\DeclareMathOperator{\sign}{sign}
\DeclareMathOperator{\diag}{diag}
\DeclareMathOperator*{\argmin}{arg\,min}

\newenvironment{keywords}{\begin{quote}\small\noindent\textbf{Keywords:}\ \ignorespaces}{\end{quote}}

\title{Sparse optimal control for a class of infinite-dimensional bilinear systems: optimality conditions and uniqueness}
\author[1]{Redouane El Mezegueldy}
\author[2]{Zakarya Dardour}
\affil[1]{MASI Laboratory, \'Ecole Normale Sup\'erieure, Sidi Mohamed Ben Abdellah University, Fez, Morocco}
\affil[2]{L2MASI, Faculty of Sciences Dhar El Mahraz, Sidi Mohamed Ben Abdellah University, Fez, Morocco}
\date{}

\begin{document}

\maketitle
\begingroup
\renewcommand{\thefootnote}{}
\footnotetext{\raggedright\emph{E-mail addresses:} \mbox{\href{mailto:redouane.elmezegueldy@usmba.ac.ma}{\nolinkurl{redouane.elmezegueldy@usmba.ac.ma}}} (R.~El Mezegueldy, corresponding author), \mbox{\href{mailto:zakarya.dardour@usmba.ac.ma}{\nolinkurl{zakarya.dardour@usmba.ac.ma}}} (Z.~Dardour).}
\endgroup

\begin{abstract}
We study optimal control of an evolution equation in which a scalar control multiplies a bounded operator acting on the state, with free dynamics given by a $C_0$-semigroup. The cost combines tracking terms, a quadratic control cost and an $L^1$ term promoting sparsity in time. We prove existence, derive optimality conditions, and show that optimal controls are thresholded functions of a switching function whose derivative is governed by the commutator of the drift and control operators. Consequently, commuting problems with a terminal cost admit no temporally sparse controls. The $L^1$ term also matters for uniqueness: two distinct controls satisfying the optimality conditions differ on a set of measure bounded from below, contained in the union of their supports, which yields, for the first time, uniqueness for every horizon in stable problems steered to zero. Second-order sufficient conditions and experiments on a bilinear heat equation complete the study.
\end{abstract}


\section{Introduction}\label{sec:intro}

In a bilinear control system the control multiplies the state instead of entering as an additive input. Such systems arise whenever an actuator changes a coefficient of the dynamics: a reaction rate in a diffusion process, the kill rate of a drug in a pharmacokinetic model, the strength of a dipole coupling in quantum mechanics, or the drift of a Fokker--Planck equation \citep{FleigGuglielmi2017}. For bounded control operators the reachable set has dense complement in infinite dimensions \citep{BallMarsdenSlemrod1982}, so exact steering is out of reach and optimal control is the natural way to formulate the task.

We consider the abstract bilinear system
\begin{equation}\label{eq:state}
y'(t)=Ay(t)+u(t)\,By(t),\quad t\in(0,T),\qquad y(0)=y_0,
\end{equation}
where $A$ generates a $C_0$-semigroup on a Hilbert space $H$, $B$ is a bounded operator on $H$, and the scalar control $u$ takes values in an interval $[a,b]$ containing $0$. Writing $y_u$ for the state associated with $u$, we minimize
\begin{equation}\label{eq:cost}
J(u)=\frac{\gamma_1}{2}\norm{y_u(T)-y_T}^2+\frac{\gamma_2}{2}\int_0^T\norm{y_u(t)-y_d(t)}^2\dd t+\frac{\alpha}{2}\norm{u}_{L^2(0,T)}^2+\beta\norm{u}_{L^1(0,T)}
\end{equation}
over the admissible set $\Uad$ of controls in $L^2(0,T)$ with values in $[a,b]$, and call this problem
\begin{equation}\label{eq:problem}
(\mathrm{P}_\beta)\qquad \min_{u\in\Uad}J(u).
\end{equation}
For $\beta=0$ this is the classical problem with a quadratic control cost. For $\beta>0$ the $L^1$ term charges the controller for every instant at which it acts, and optimal controls switch off on whole time intervals. This is the desired behaviour when actuation is expensive or harmful, as with drug doses or radio-frequency pulses, and it raises two questions: when does sparsity in time occur, and how does the $L^1$ term change the properties of the problem?

Sparse controls are preferred in engineering, but they make the mathematical problem harder. The $L^1$ norm is not differentiable, and it fails to be differentiable precisely at the controls of interest, those that vanish on whole intervals. The first-order condition then involves the subdifferential of the $L^1$ norm, the pointwise control law acquires kinks at the thresholds $\pm\beta$, second-order conditions must take the kink of the cost at $u=0$ into account, and Newton-type methods have to be replaced by nonsmooth ones. It is therefore not clear in advance which properties of the smooth problem survive.

We show that the classical theory carries over, and that the nonsmooth term also brings something that the smooth cost does not provide. The uniqueness results available for bilinear problems are local in time: they require the horizon $T$ to be small relative to the initial and desired data and to the regularization parameter. With the $L^1$ term, an exponentially stable system steered to zero has a unique optimal control for every horizon, provided $\alpha$ exceeds a threshold that depends neither on the horizon nor on the control bounds. For the smooth cost, our argument gives such a result only when the control bounds are small compared with the decay rate of the semigroup. This gain is somewhat surprising, because the $L^1$ term does not convexify the problem (Example~\ref{ex:nonunique}). It comes instead from a bound on the $L^1$ norm of optimal controls that does not depend on the horizon, combined with the decay of state and adjoint.

Optimal control of bilinear distributed systems is well developed. The classical results give existence of optimal controls, a first-order optimality system, and uniqueness under a smallness condition on the final time \citep{BradleyLenhart1997,BouhamedElKaboussBouzahir2022} or on the initial state \citep{AddouBenbrik2002}; bounded and unbounded control sets are treated by \citet{ZerrikAitAadiOuhafsa2022}. When the cost is linear in the control, \citet{AronnaBonnansKroner2018} derived second-order conditions for semigroup models through the Goh transform, and \citet{CasasWachsmuthWachsmuth2018} analysed bang-bang bilinear control of elliptic equations.

Sparse optimal control was introduced in PDE-constrained optimization by \citet{Stadler2009} and has since been developed in many directions \citep{WachsmuthWachsmuth2011,CasasHerzogWachsmuth2012,HerzogStadlerWachsmuth2012,CasasHerzogWachsmuth2017,HoppeNeitzel2022}. Closest to our setting are results on sparsity in time: for the stabilization of semilinear parabolic equations, an $L^1$-type cost makes the optimal control vanish after a finite time \citep{CasasKunisch2017}, and infinite-horizon sparse control of finite-dimensional systems is studied by \citet{KaliseKunischRao2017}. In control engineering, maximum hands-off control computes the sparsest control in time by $L^1$ minimization \citep{NagaharaQuevedoNesic2016,ChatterjeeNagaharaQuevedoRao2016}, and $L^1$-minimal controls of control-affine systems are analysed by \citet{VossenMaurer2006}. For bilinear systems, $L^1$ costs have mostly been used in quantum control, where they produce sparse pulses \citep{CiaramellaBorzi2016} or fields supported on few frequencies \citep{FrieseckeHennekeKunisch2018}, and where the $L^1$ norm of the field measures its effect on the state \citep{BoussaidCaponigroChambrion}. In biomedical models, objectives that are linear in the control are argued to describe drug cost and toxicity better than quadratic ones \citep{LedzewiczSchattler2020}.

Two further lines of work are relevant here. The structure of an optimal control is read off its switching function $\sigma=\ip{By}{p}$, where $p$ is the adjoint state; for control-affine ordinary differential equations, the derivative of the switching function is given by a Lie bracket of the drift and the control vector field \citep{SchattlerLedzewicz2012}. Long horizons with zero targets lead to the turnpike phenomenon, in which optimal controls for different horizons nearly coincide away from the final time \citep{TrelatZuazua2015,GruneSchallerSchiela2020}; infinite-horizon bilinear problems are treated by \citet{BreitenKunischPfeiffer2018}.

To our knowledge, the abstract bilinear problem with a combined $L^2$--$L^1$ cost has not been analysed, and the effect of the $L^1$ term on uniqueness has not been examined. We state every result for $\beta\ge0$, so that each theorem shows what already holds for the smooth cost and what the $L^1$ term adds.

Our first contribution is the basic theory of the problem. We prove that optimal controls exist and derive an optimality system in which the optimal control is a projected soft-thresholding of the switching function, with explicit critical weights beyond which the zero control satisfies this optimality system or is the unique solution. Second, we establish an identity for the derivative of the switching function, the counterpart of the Lie-bracket formula, which shows that sparsity in time is driven by the commutator $BA-AB$. When $A$ and $B$ commute and only the final state is penalized, every control in $\Uad$ that satisfies the optimality system is constant in time, so no sparsity in time can occur.

Our third and main contribution concerns uniqueness. We prove a separation estimate: two distinct controls in $\Uad$ that satisfy the optimality system differ on a set of measure at least $\alpha/L_\sigma$, where $L_\sigma$ is the Lipschitz constant of the switching function with respect to the control in the $L^1$ norm, and this set lies in the union of their supports. For $\beta>0$ the supports are the sets where $\abs{\sigma}$ exceeds $\beta$, so short supports force uniqueness. For exponentially stable systems steered to zero, this yields a unique optimal control for every horizon, with constants independent of the horizon and of the control bounds. To our knowledge, this is the first uniqueness result for bilinear optimal control problems that holds on every time interval $[0,T]$, whereas the earlier results \citep{BradleyLenhart1997,AddouBenbrik2002,BouhamedElKaboussBouzahir2022} are local in time.

Fourth, we give second-order sufficient conditions together with an a posteriori test in which the measure of the support replaces the length of the horizon. Finally, we propose a positivity-preserving discretization with an exact discrete adjoint and report experiments on a bilinear heat equation that confirm these properties, including mesh independence and a direct check of the second-order condition.

The paper is organized as follows. Section~\ref{sec:setting} states the standing assumptions and collects the properties of the state equation used throughout: well-posedness, Lipschitz dependence of the state on the control in the $L^1$ norm, continuity from weak to strong convergence, differentiability, and a duality identity between state and adjoint. Section~\ref{sec:existence} proves the existence of optimal controls. Section~\ref{sec:os} derives the optimality system, the soft-thresholding law and the critical weights for the zero control. Section~\ref{sec:switching} establishes the switching identity and treats the commuting case. Section~\ref{sec:uniqueness} contains the separation estimate and the uniqueness results, including uniqueness for every horizon in stable problems steered to zero, together with a scalar example of non-uniqueness. Section~\ref{sec:ssc} gives second-order sufficient conditions and the a posteriori test. Section~\ref{sec:numerics} reports numerical experiments on a bilinear heat equation, and Section~\ref{sec:conclusion} concludes.

\section{Setting and preliminaries}\label{sec:setting}

The analysis rests on four properties of the control-to-state map: an a priori bound, Lipschitz dependence on the control in the $L^1$ norm, continuity from weak to strong convergence, and smoothness. A duality identity between state and adjoint completes the toolbox. These facts are classical, and most of them go back to the analysis of distributed bilinear systems by \citet{BallMarsdenSlemrod1982}; we state them in the form used later, refer to the literature for their proofs, and give arguments only where explicit constants are needed in Sections~\ref{sec:uniqueness} and~\ref{sec:ssc}. For background on semigroups we refer to \citet{EngelNagel2000} and \citet{Pazy1983}, and for optimal control of evolution equations to \citet{LiYong1995} and \citet{Troltzsch2010}.

Throughout the paper, $H$ is a real Hilbert space with inner product $\ip{\cdot}{\cdot}$ and norm $\norm{\cdot}$, and the operator $A\colon D(A)\subset H\to H$ generates a $C_0$-semigroup $(S(t))_{t\ge0}$ on $H$. We set $M_T:=\sup_{0\le t\le T}\norm{S(t)}$ and recall that the adjoint operators $S(t)^*$ form a $C_0$-semigroup generated by $A^*$ \citep[Section~1.10]{Pazy1983}. The control operator satisfies $B\in\mathcal{L}(H)$, the space of bounded linear operators on $H$. The admissible set is $\Uad:=\{u\in L^2(0,T)\colon a\le u(t)\le b\ \text{for a.e. } t\}$, where $a\le0\le b$ and $a<b$, and we set $\rho:=\max\{-a,b\}$. The data satisfy $y_0,y_T\in H$ and $y_d\in L^2(0,T;H)$, and the weights satisfy $\gamma_1,\gamma_2\ge0$, $\gamma_1+\gamma_2>0$, $\alpha>0$ and $\beta\ge0$. These standing assumptions are used in all results without further mention; additional hypotheses are stated where they are needed.

A mild solution of the state equation \eqref{eq:state} is a function $y\in C([0,T];H)$ such that
\begin{equation}\label{eq:mild}
y(t)=S(t)y_0+\int_0^t S(t-s)\,u(s)\,By(s)\dd s,\qquad t\in[0,T].
\end{equation}
We write $J=F+\beta j$, where $j(u):=\norm{u}_{L^1(0,T)}$ and $F$ collects the other three terms. The norm of $C([0,T];H)$ or $C([0,T])$ is $\norm{\cdot}_C$, and $\abs{E}$ is the Lebesgue measure of a set $E\subset[0,T]$. The constants
\begin{equation}\label{eq:constants}
E_T:=e^{M_T\norm{B}\rho T},\qquad K_y:=M_T\norm{y_0}E_T,\qquad L_y:=M_T\norm{B}K_yE_T
\end{equation}
appear throughout.

We first show that the state equation is well posed for every integrable control.

\begin{lemma}\label{lem:wellposed}
For every $u\in L^1(0,T)$, \eqref{eq:state} has a unique mild solution $y_u\in C([0,T];H)$, and
\[
\norm{y_u(t)}\le M_T\norm{y_0}\exp\Bigl(M_T\norm{B}\int_0^t\abs{u(s)}\dd s\Bigr),\qquad t\in[0,T].
\]
In particular $\norm{y_u}_C\le K_y$ for all $u\in\Uad$.
\end{lemma}

\begin{proof}
For $u\in L^1(0,T)$, existence and uniqueness of the mild solution in $C([0,T];H)$ is a special case of \citet[Theorem~2.5]{BallMarsdenSlemrod1982}, which allows controls in $L^1(0,T)$ and covers every bounded linear operator $B$. For bounded controls, such as those in $\Uad$, it is also a special case of the existence theory for semilinear evolution equations \citep[Chapter~2, Proposition~5.3]{LiYong1995}. The bound follows from $\norm{y_u(t)}\le M_T\norm{y_0}+M_T\norm{B}\int_0^t\abs{u(s)}\norm{y_u(s)}\dd s$ and Gronwall's lemma.
\end{proof}

The state depends Lipschitz continuously on the control, measured in the $L^1$ norm. This is a quantitative form of the continuity of $u\mapsto y_u$ from $L^1(0,T)$ into $C([0,T];H)$ contained in \citet[Theorem~2.5]{BallMarsdenSlemrod1982}; we give the short proof because the constant $L_y$ is used in Section~\ref{sec:uniqueness}.

\begin{lemma}\label{lem:lipschitz}
For $u_1,u_2\in\Uad$,
\[
\norm{y_{u_1}-y_{u_2}}_C\le L_y\norm{u_1-u_2}_{L^1(0,T)}.
\]
\end{lemma}

\begin{proof}
The difference $\delta y:=y_{u_1}-y_{u_2}$ satisfies $\delta y(t)=\int_0^tS(t-s)\bigl[u_1(s)B\delta y(s)+(u_1-u_2)(s)By_{u_2}(s)\bigr]\dd s$. Therefore $\norm{\delta y(t)}\le M_T\norm{B}K_y\norm{u_1-u_2}_{L^1}+M_T\norm{B}\int_0^t\abs{u_1}\norm{\delta y}\dd s$, and Gronwall's lemma with $\norm{u_1}_{L^1}\le\rho T$ gives the claim.
\end{proof}

The control enters this estimate only through its $L^1$ norm. This is what later allows the measure of a support, rather than the length of the horizon, to appear in the uniqueness results.

Continuity from weakly convergent controls to uniformly convergent states is due to \citet[Theorem~3.6]{BallMarsdenSlemrod1982}, where it is the key step towards their noncontrollability result. It rests on the following convergence lemma, which we also use in Section~\ref{sec:ssc} for the linearized equation.

\begin{lemma}[{\citealp[Lemma~3.7]{BallMarsdenSlemrod1982}}]\label{lem:bms}
Let $u\in L^1(0,T)$ and let $w_n\rightharpoonup0$ weakly in $L^1(0,T)$. Then
\[
\lim_{n\to\infty}\,\sup_{0\le t\le T}\Bigl\|\int_0^tS(t-s)\,w_n(s)\,By_u(s)\dd s\Bigr\|=0.
\]
\end{lemma}

The proof in \citet{BallMarsdenSlemrod1982} uses only the strong continuity of the semigroup, the continuity of $s\mapsto By_u(s)$ and properties of weakly convergent sequences in $L^1(0,T)$. Together with Gronwall's lemma, it gives the continuity of the control-to-state map from weak to strong convergence.

\begin{lemma}[{\citealp[Theorem~3.6]{BallMarsdenSlemrod1982}}]\label{lem:weakstrong}
If $u_n\rightharpoonup u$ weakly in $L^1(0,T)$, in particular if $u_n\rightharpoonup u$ weakly in $L^2(0,T)$, then $y_{u_n}\to y_u$ in $C([0,T];H)$.
\end{lemma}

\begin{proof}
Write $y_n:=y_{u_n}$ and $y:=y_u$. Subtracting the mild formulations \eqref{eq:mild} gives
\[
y_n(t)-y(t)=\int_0^tS(t-s)\,u_n(s)\,B\bigl(y_n(s)-y(s)\bigr)\dd s+\int_0^tS(t-s)\,(u_n-u)(s)\,By(s)\dd s.
\]
Let $\varepsilon_n$ be the supremum over $t\in[0,T]$ of the norm of the last integral; $\varepsilon_n\to0$ by Lemma~\ref{lem:bms}. Weakly convergent sequences are bounded, so $c:=\sup_n\norm{u_n}_{L^1}<\infty$, and Gronwall's lemma gives $\norm{y_n-y}_C\le\varepsilon_ne^{M_T\norm{B}c}\to0$.
\end{proof}

The control-to-state map is moreover smooth, and its derivatives solve linearized equations.

\begin{lemma}\label{lem:diff}
The map $G\colon L^1(0,T)\to C([0,T];H)$, $G(u):=y_u$, is of class $C^\infty$; in particular, its restriction to $L^2(0,T)$ is twice continuously Fr\'echet differentiable. For $u,v\in L^1(0,T)$, $z_v:=G'(u)v$ and $w_v:=G''(u)[v,v]$ are the mild solutions of
\begin{equation}\label{eq:linearized}
z'=Az+u\,Bz+v\,By_u,\quad z(0)=0;\qquad w'=Aw+u\,Bw+2v\,Bz_v,\quad w(0)=0.
\end{equation}
Moreover, with $\kappa:=M_T\norm{B}$,
\begin{equation}\label{eq:zbound}
\norm{z_v}_C\le\kappa\,\norm{y_u}_C\,e^{\kappa\norm{u}_{L^1(0,T)}}\,\norm{v}_{L^1(0,T)},
\end{equation}
and in particular $\norm{z_v}_C\le L_y\norm{v}_{L^1(0,T)}$ for $u\in\Uad$.
\end{lemma}

\begin{proof}
The smoothness of $G$ is part of \citet[Theorem~2.5]{BallMarsdenSlemrod1982}, and the characterization of $G'(u)v$ through the linearized equation is their Corollary~2.4; \citet[Section~2]{AronnaBonnansKroner2018} prove that $G$ is of class $C^\infty$ from $L^1(0,T)$ into $C([0,T];H)$ by the implicit function theorem. With $V(u,y)(t):=\int_0^tS(t-s)\,u(s)\,By(s)\dd s$, which is bilinear and continuous from $L^1(0,T)\times C([0,T];H)$ into $C([0,T];H)$, the state satisfies $G(u)=S(\cdot)y_0+V(u,G(u))$. Differentiating this identity once and twice in the direction $v$ gives $z_v=V(u,z_v)+V(v,y_u)$ and $w_v=V(u,w_v)+2V(v,z_v)$, which are the mild forms of \eqref{eq:linearized}. Finally, $\norm{z_v(t)}\le\kappa\norm{y_u}_C\norm{v}_{L^1}+\kappa\int_0^t\abs{u}\norm{z_v}\dd s$, and Gronwall's lemma gives \eqref{eq:zbound}. For $u\in\Uad$, $\norm{y_u}_C\le K_y$ and $e^{\kappa\norm{u}_{L^1}}\le E_T$.
\end{proof}

Finally, a duality identity links the state and adjoint equations; it is used below to compute derivatives of the cost and of the switching function. It is the integration-by-parts formula for mild solutions and can also be derived from the equivalence of weak and mild solutions \citep{Ball1977}; see \citet[Section~2]{AronnaBonnansKroner2018} for such identities in the bilinear setting. We include the short proof.

\begin{lemma}\label{lem:duality}
Let $0\le t_0<T$, $f,g\in L^1(t_0,T;H)$ and $z_0,q_T\in H$. For $t\in[t_0,T]$ let $z(t):=S(t-t_0)z_0+\int_{t_0}^tS(t-s)f(s)\dd s$ and $q(t):=S(T-t)^*q_T+\int_t^TS(s-t)^*g(s)\dd s$. Then
\begin{equation}\label{eq:duality}
\ip{z(T)}{q_T}-\ip{z(t_0)}{q(t_0)}=\int_{t_0}^T\bigl(\ip{f(s)}{q(s)}-\ip{z(s)}{g(s)}\bigr)\dd s.
\end{equation}
\end{lemma}

\begin{proof}
Insert the formula for $q$ into $\int\ip{f}{q}\dd s$ and apply Fubini's theorem:
\[
\int_{t_0}^T\ip{f}{q}\dd s=\Bigl\langle\int_{t_0}^TS(T-s)f(s)\dd s,\;q_T\Bigr\rangle+\int_{t_0}^T\Bigl\langle\int_{t_0}^rS(r-s)f(s)\dd s,\;g(r)\Bigr\rangle\dd r.
\]
The inner integrals equal $z(T)-S(T-t_0)z_0$ and $z(r)-S(r-t_0)z_0$. The terms containing $z_0$ add up to $\ip{z_0}{S(T-t_0)^*q_T+\int_{t_0}^TS(r-t_0)^*g(r)\dd r}=\ip{z(t_0)}{q(t_0)}$, which gives \eqref{eq:duality}.
\end{proof}

\section{Existence of an optimal control}\label{sec:existence}

The direct method, combined with Lemma~\ref{lem:weakstrong}, shows that $(\mathrm{P}_\beta)$ is solvable for every $\beta\ge0$.

\begin{theorem}\label{thm:existence}
Problem $(\mathrm{P}_\beta)$ has at least one solution $\bar u\in\Uad$ for every $\beta\ge0$. The same holds for $\alpha=0$.
\end{theorem}

\begin{proof}
$J$ is bounded below by $0$. Let $(u_n)\subset\Uad$ be a minimizing sequence. Since $\Uad$ is bounded, closed and convex in $L^2(0,T)$, a subsequence, not relabeled, converges weakly to some $\bar u\in\Uad$. By Lemma~\ref{lem:weakstrong}, $y_{u_n}\to y_{\bar u}$ in $C([0,T];H)$, so both tracking terms of \eqref{eq:cost} converge. The terms $\frac{\alpha}{2}\norm{u}^2$ and $\beta\norm{u}_{L^1}$ are convex and continuous on $L^2(0,T)$, hence weakly lower semicontinuous. Therefore $J(\bar u)\le\liminf J(u_n)=\inf J$, and $\bar u$ solves $(\mathrm{P}_\beta)$.
\end{proof}

\begin{remark}\label{rem:unconstrained}
If $\Uad$ is replaced by $L^2(0,T)$, existence still holds for $\alpha>0$, because $J(u)\ge\frac{\alpha}{2}\norm{u}^2$ keeps minimizing sequences bounded. We keep the box constraints because they make the constants \eqref{eq:constants} uniform over all admissible controls. Results whose constants do not involve $a$ and $b$, such as Theorem~\ref{thm:longhorizon}, hold for $\Uad=L^2(0,T)$ as well, with $P$ replaced by the identity in \eqref{eq:law}.
\end{remark}

\begin{remark}\label{rem:nocompact}
As observed by \citet{BallMarsdenSlemrod1982}, Lemma~\ref{lem:weakstrong} uses neither compactness of $S(t)$ nor compactness of $B$; it only uses that the control is a scalar function of time. Theorem~\ref{thm:existence} therefore covers parabolic, hyperbolic and Schr\"odinger-type equations alike.
\end{remark}

\section{Optimality system and sparsity structure}\label{sec:os}

Every local solution satisfies a pointwise soft-thresholding law driven by the switching function $\sigma=\ip{By}{p}$. The law has a dead zone exactly when $\beta>0$, and it is this dead zone that produces sparsity.

For $u\in\Uad$, the adjoint state $p_u\in C([0,T];H)$ is the mild solution of
\begin{equation}\label{eq:adjoint}
-p'(t)=A^*p(t)+u(t)\,B^*p(t)+\gamma_2\bigl(y_u(t)-y_d(t)\bigr),\qquad p(T)=\gamma_1\bigl(y_u(T)-y_T\bigr),
\end{equation}
that is, $p_u(t)=S(T-t)^*\gamma_1(y_u(T)-y_T)+\int_t^TS(s-t)^*\bigl[u(s)B^*p_u(s)+\gamma_2(y_u(s)-y_d(s))\bigr]\dd s$. The adjoint equation is linear in $p$, with the integrable coefficient $u$ and an integrable forcing term, so existence and uniqueness follow as for the state equation, backward in time, and Gronwall's lemma gives, for all $u\in\Uad$,
\[
\norm{p_u}_C\le K_p:=M_TE_T\Bigl(\gamma_1\bigl(K_y+\norm{y_T}\bigr)+\gamma_2\bigl(TK_y+\sqrt T\norm{y_d}_{L^2(0,T;H)}\bigr)\Bigr).
\]
The switching function is
\begin{equation}\label{eq:sigma}
\sigma_u(t):=\ip{By_u(t)}{p_u(t)},\qquad t\in[0,T],
\end{equation}
and $\sigma_u\in C([0,T])$ with $\norm{\sigma_u}_C\le\norm{B}K_yK_p$.

The switching function is the $L^2$ gradient of the tracking part of the cost.

\begin{proposition}\label{prop:derivative}
$F$ is twice continuously Fr\'echet differentiable on $L^2(0,T)$, and for $u,v\in L^2(0,T)$
\begin{equation}\label{eq:derivative}
F'(u)v=\int_0^T\bigl(\alpha u(t)+\sigma_u(t)\bigr)v(t)\dd t.
\end{equation}
\end{proposition}

\begin{proof}
The computation is the standard adjoint calculus, carried out for bilinear systems with a cost linear in the control in \citet[Section~2]{AronnaBonnansKroner2018}; we recall it to fix the notation. Differentiability follows from Lemma~\ref{lem:diff} and the chain rule. With $z_v$ from \eqref{eq:linearized}, the derivative of the tracking terms in the direction $v$ is $\gamma_1\ip{y_u(T)-y_T}{z_v(T)}+\gamma_2\int_0^T\ip{y_u-y_d}{z_v}\dd t$. Apply Lemma~\ref{lem:duality} with $t_0=0$, $z_0=0$, $f=uBz_v+vBy_u$, $q_T=\gamma_1(y_u(T)-y_T)$ and $g=uB^*p_u+\gamma_2(y_u-y_d)$, so that $q=p_u$. The terms containing $u$ cancel because $\ip{uBz_v}{p_u}=\ip{z_v}{uB^*p_u}$, and \eqref{eq:duality} gives
\[
\gamma_1\ip{y_u(T)-y_T}{z_v(T)}+\gamma_2\int_0^T\ip{y_u-y_d}{z_v}\dd t=\int_0^Tv\,\ip{By_u}{p_u}\dd t.
\]
Adding the derivative of $\frac{\alpha}{2}\norm{u}^2$ gives \eqref{eq:derivative}.
\end{proof}

We use the following notation. For $s\in\R$ let $S_\beta(s):=\sign(s)\max\{\abs{s}-\beta,0\}$, the soft-thresholding at level $\beta$, and let $P$ denote the projection of $\R$ onto $[a,b]$. For a continuous function $u$ on $[0,T]$ we write $\supp u:=\{t\in[0,T]\colon u(t)\neq0\}$; this relatively open set, not its closure, is what we call the support, and $\abs{\supp u}$ is the time during which $u$ acts.

We can now state the optimality system.

\begin{theorem}\label{thm:os}
Let $\bar u$ be a local solution of $(\mathrm{P}_\beta)$ in $L^2(0,T)$, that is, $J(\bar u)\le J(u)$ for some $\varepsilon>0$ and all $u\in\Uad$ with $\norm{u-\bar u}_{L^2}\le\varepsilon$. Let $\bar y$, $\bar p$ and $\bar\sigma$ be its state, adjoint and switching function. Then, after choosing the continuous representative,
\begin{equation}\label{eq:law}
\bar u(t)=P\Bigl(-\frac{1}{\alpha}S_\beta\bigl(\bar\sigma(t)\bigr)\Bigr)\qquad\text{for all } t\in[0,T].
\end{equation}
In what follows, a control $\bar u\in\Uad$ is said to satisfy \eqref{eq:law} if \eqref{eq:law} holds with $\bar\sigma=\sigma_{\bar u}$; by the theorem, every local solution does.
\end{theorem}

\begin{proof}
The derivation is the one known for elliptic problems with $L^1$ cost \citep{Stadler2009,CasasHerzogWachsmuth2012}, so we only recall it. Let $u\in\Uad$ and $u_s:=\bar u+s(u-\bar u)$ for $s\in(0,1]$. Local optimality and convexity of $j$ give, for small $s$, $0\le J(u_s)-J(\bar u)\le F(u_s)-F(\bar u)+s\beta\bigl(j(u)-j(\bar u)\bigr)$. Dividing by $s$ and letting $s\downarrow0$ yields, by \eqref{eq:derivative},
\begin{equation}\label{eq:vi}
\int_0^T\bigl(\alpha\bar u+\bar\sigma\bigr)(u-\bar u)\dd t+\beta\int_0^T\bigl(\abs{u}-\abs{\bar u}\bigr)\dd t\ge0\qquad\text{for all }u\in\Uad.
\end{equation}
Since $\Uad$ is defined by pointwise bounds, \eqref{eq:vi} localizes in the usual way (take $u=w$ on a measurable set and $u=\bar u$ elsewhere, with $w\in[a,b]$): for a.e.\ $t$, $(\alpha\bar u(t)+\bar\sigma(t))(w-\bar u(t))+\beta(\abs{w}-\abs{\bar u(t)})\ge0$ for all $w\in[a,b]$. By the gradient inequality for $\frac{\alpha}{2}w^2$, $\bar u(t)$ minimizes the strictly convex function $\varphi_t(w):=\frac{\alpha}{2}w^2+\bar\sigma(t)w+\beta\abs{w}$ over $[a,b]$, so it is the projection onto $[a,b]$ of the unconstrained minimizer $-S_\beta(\bar\sigma(t))/\alpha$. This gives \eqref{eq:law} for a.e.\ $t$; it is the soft-thresholding formula of \citet[Section~3]{Stadler2009}, with the switching function in place of the adjoint state and followed by the projection onto $[a,b]$. The right-hand side of \eqref{eq:law} is continuous in $t$, since $\bar\sigma$ is continuous and $s\mapsto P(-S_\beta(s)/\alpha)$ is Lipschitz; this yields the continuous representative.
\end{proof}

The pointwise law has several immediate consequences.

\begin{corollary}\label{cor:structure}
Let $\bar u\in\Uad$ satisfy \eqref{eq:law}, and set $\Phi(s):=P(-S_\beta(s)/\alpha)$ for $s\in\R$. Then the following hold.
\begin{enumerate}[label=\textup{(\arabic*)},leftmargin=*]
\item $\bar u\in C([0,T])$ and $\abs{\bar u(t)}\le\max\{\abs{\bar\sigma(t)}-\beta,0\}/\alpha$ for all $t$;
\item $\abs{\Phi(s_1)-\Phi(s_2)}\le\abs{s_1-s_2}/\alpha$ for all $s_1,s_2\in\R$;
\item if $a<0<b$, then $\bar u(t)=0$ if and only if $-\beta\le\bar\sigma(t)\le\beta$; for nonnegative controls ($a=0$), $\bar u(t)=0$ if and only if $\bar\sigma(t)\ge-\beta$;
\item if $a<0<b$ and $\beta=0$, then $\{\bar u=0\}=\{\bar\sigma=0\}$: the smooth cost switches the control off on an interval only where $\bar\sigma$ vanishes identically. If $\beta>0$, then $\supp\bar u=\{\abs{\bar\sigma}>\beta\}$, and $\bar u$ is zero wherever the switching function is below the price $\beta$.
\end{enumerate}
\end{corollary}

\begin{proof}
Items (1) and (2) hold because $S_\beta$ and $P$ are $1$-Lipschitz \citep{BauschkeCombettes2017} and $P(0)=0$, with $\abs{P(x)}\le\abs{x}$. For item (3), $P(x)=0$ if and only if $x=0$ when $a<0<b$, and $S_\beta(s)=0$ if and only if $\abs{s}\le\beta$; the case $a=0$ is checked in the same way. Item (4) follows from item (3) and the continuity of $\bar\sigma$.
\end{proof}

In particular, the zero control satisfies \eqref{eq:law} exactly when $\beta$ dominates its switching function, and it is the only control in $\Uad$ satisfying \eqref{eq:law} when $\beta$ dominates all admissible switching functions.

\begin{proposition}\label{prop:critical}
Assume $a<0<b$. Let $\sigma_0$ be the switching function of the zero control, computed from the free state $S(t)y_0$ and its adjoint, and set
\[
\beta_0:=\norm{\sigma_0}_C,\qquad \bar\beta:=\sup_{u\in\Uad}\norm{\sigma_u}_C\le\norm{B}K_yK_p.
\]
Then $u\equiv0$ satisfies \eqref{eq:law} if and only if $\beta\ge\beta_0$. If $\beta\ge\bar\beta$, the zero control is the only control in $\Uad$ satisfying \eqref{eq:law} and therefore the unique solution of $(\mathrm{P}_\beta)$.
\end{proposition}

\begin{proof}
By Corollary~\ref{cor:structure}(3), $u\equiv0$ satisfies \eqref{eq:law} if and only if $\abs{\sigma_0(t)}\le\beta$ for all $t$. If $\beta\ge\bar\beta$ and $\bar u\in\Uad$ satisfies \eqref{eq:law}, then $\abs{\bar\sigma(t)}\le\bar\beta\le\beta$ for all $t$, hence $\bar u\equiv0$. Solutions exist by Theorem~\ref{thm:existence} and satisfy \eqref{eq:law} by Theorem~\ref{thm:os}.
\end{proof}

For $\beta=0$ the zero control satisfies \eqref{eq:law} only if $\sigma_0\equiv0$; the smooth cost has no threshold. For convex linear-quadratic elliptic problems, \citet[Lemma~3.1]{Stadler2009} shows that the zero control is the solution as soon as $\beta$ is at least the supremum norm of the adjoint state of the zero control, the analogue of $\beta_0$. Without convexity, $\beta\ge\beta_0$ only ensures that the zero control satisfies \eqref{eq:law}, and Proposition~\ref{prop:critical} guarantees its optimality only for the larger weight $\bar\beta$. Between $\beta_0$ and $\bar\beta$, the zero control can coexist with other controls satisfying \eqref{eq:law} (Example~\ref{ex:nonunique}). For nonnegative controls ($a=0$), Corollary~\ref{cor:structure}(3) shows in the same way that $u\equiv0$ satisfies \eqref{eq:law} if and only if $\beta\ge\max_{t\in[0,T]}\max\{-\sigma_0(t),0\}$.

For global solutions, increasing $\beta$ trades tracking accuracy for a smaller control.

\begin{proposition}\label{prop:monotone}
Let $0\le\beta_1<\beta_2$ and let $u_i$ be a global solution of $(\mathrm{P}_{\beta_i})$, $i=1,2$. Then $\norm{u_2}_{L^1}\le\norm{u_1}_{L^1}$ and $F(u_1)\le F(u_2)$.
\end{proposition}

\begin{proof}
Optimality gives $F(u_1)+\beta_1\norm{u_1}_{L^1}\le F(u_2)+\beta_1\norm{u_2}_{L^1}$ and $F(u_2)+\beta_2\norm{u_2}_{L^1}\le F(u_1)+\beta_2\norm{u_1}_{L^1}$. Adding them gives $(\beta_2-\beta_1)(\norm{u_2}_{L^1}-\norm{u_1}_{L^1})\le0$. Inserting this into the first inequality gives $F(u_1)\le F(u_2)$.
\end{proof}

This is the classical monotonicity argument for Tikhonov regularization with respect to the regularization parameter \citep{TikhonovLeonovYagola1998}. No convexity is used: the $L^1$ norm of global solutions decreases and the tracking part increases with $\beta$, even when \eqref{eq:law} has several solutions in $\Uad$.

\begin{remark}\label{rem:general}
More generally, $\frac{\alpha}{2}w^2+\beta\abs{w}$ can be replaced by any $\alpha$-strongly convex function $g$ with $g(0)=0$. The proof of Theorem~\ref{thm:os} then gives $\bar u(t)=\Phi_g(\bar\sigma(t))$ with $\Phi_g(s):=\argmin_{w\in[a,b]}\{g(w)+sw\}$, and $\Phi_g$ is $(1/\alpha)$-Lipschitz. When $a<0<b$, $\Phi_g(s)=0$ if and only if $-s\in\partial g(0)$, so $\Phi_g$ has a dead zone exactly when $\partial g(0)$ is a nondegenerate interval, that is, when $g$ has a kink at $0$. The results below that use only the Lipschitz continuity of $\Phi_g$ therefore hold for smooth costs too, while those that use the dead zone need the kink.
\end{remark}

\section{Structure of the switching function}\label{sec:switching}

This section explains when temporal sparsity can occur at all. If $B$ leaves $D(A)$ invariant, the switching function is differentiable and its derivative is governed by the commutator $BA-AB$. In particular, if $A$ and $B$ commute and only the final state is penalized, $\sigma$ is constant, and no strongly convex control cost, smooth or sparse, can produce zero intervals.

Throughout this section we assume, in addition, that $y_0\in D(A)$ and $B(D(A))\subset D(A)$; the condition on $y_0$ is not needed in Corollary~\ref{cor:commuting} and Example~\ref{ex:identity}. By the closed graph theorem, $B$ is then bounded on $D(A)$ with the graph norm, and the commutator
\begin{equation}\label{eq:commutator}
C:=BA-AB\colon D(A)\to H
\end{equation}
is bounded from $D(A)$, with the graph norm, into $H$.

The switching function is then differentiable, and its derivative is expressed through $C$. Formally, differentiating $\sigma=\ip{By}{p}$ along \eqref{eq:state} and \eqref{eq:adjoint} gives
\[
\sigma'=\ip{BAy+uB^2y}{p}-\ip{By}{A^*p+uB^*p+\gamma_2(y-y_d)}=\ip{Cy}{p}-\gamma_2\ip{By}{y-y_d}.
\]
The control drops out because $B$ commutes with itself, and the drift survives only through the commutator. This computation is only formal, because $p(t)$ need not belong to $D(A^*)$. The proof below therefore keeps $A$ on the side of the state, where the assumptions keep $By$ in $D(A)$, and replaces the product rule by the duality identity of Lemma~\ref{lem:duality}.

\begin{theorem}\label{thm:switching}
For every $u\in\Uad$, $\sigma_u\in H^1(0,T)$ and, for a.e.\ $t\in(0,T)$,
\begin{equation}\label{eq:switching}
\sigma_u'(t)=\ip{Cy_u(t)}{p_u(t)}-\gamma_2\ip{By_u(t)}{y_u(t)-y_d(t)}.
\end{equation}
If moreover $y_d\in C([0,T];H)$, then $\sigma_u\in C^1([0,T])$.
\end{theorem}

\begin{proof}
Write $y:=y_u$ and $p:=p_u$.

\emph{Step 1: $y$ is a strong solution.} The part of $S(t)$ in $X_1:=D(A)$, with the graph norm, is a $C_0$-semigroup on $X_1$ \citep[Section~II.5]{EngelNagel2000}, and $B$ restricts to a bounded operator on $X_1$. Lemma~\ref{lem:wellposed} in $X_1$ gives a mild solution in $C([0,T];X_1)$, which coincides with $y$ by uniqueness in $H$. For $z\in D(A)$ we have $S(\tau)z=z+\int_0^\tau S(r)Az\dd r$ \citep[Theorem~1.2.4]{Pazy1983}. Inserting this into \eqref{eq:mild}, with $z=y_0$ and $z=u(s)By(s)$, and using Fubini's theorem gives $y(t)=y_0+\int_0^t(Ay+uBy)\dd s$. Hence $y\in W^{1,1}(0,T;H)$ and $y'=Ay+uBy$ a.e.

\emph{Step 2: an equation for $By$.} Let $w:=By$. Since $y\in C([0,T];X_1)$ and $B$ is bounded on $X_1$, we have $w\in C([0,T];X_1)$; since $y\in W^{1,1}(0,T;H)$ and $B$ is bounded on $H$, we have $w\in W^{1,1}(0,T;H)$ and $w'=By'$. Hence, a.e.\ on $(0,T)$,
\[
w'=BAy+u\,B^2y=Aw+f,\qquad f:=u\,Bw+Cy,
\]
and $f\in L^1(0,T;H)$ because $Cy\in C([0,T];H)$. We claim that $w$ is the mild solution of this equation from any initial time $t_0\in[0,T)$, that is, $w(t)=S(t-t_0)w(t_0)+\int_{t_0}^tS(t-s)f(s)\dd s$ for $t\in[t_0,T]$. Indeed, for $x\in D(A^*)$ the function $s\mapsto\ip{w(s)}{x}$ is absolutely continuous and, since $w(s)\in D(A)$,
\[
\frac{\mathrm d}{\mathrm ds}\ip{w(s)}{x}=\ip{Aw(s)+f(s)}{x}=\ip{w(s)}{A^*x}+\ip{f(s)}{x}\qquad\text{for a.e. } s.
\]
Thus $w$ is a weak solution in the sense of \citet{Ball1977}, and by the main result of that paper it coincides with the mild solution.

\emph{Step 3: duality.} Apply Lemma~\ref{lem:duality} on $[t,T]$ with $z=w$, $f=uBw+Cy$ and $q=p$, that is, $q_T=\gamma_1(y(T)-y_T)$ and $g=uB^*p+\gamma_2(y-y_d)$. The terms containing $u$ cancel, and
\[
\sigma_u(T)-\sigma_u(t)=\int_t^T\Bigl(\ip{Cy}{p}-\gamma_2\ip{By}{y-y_d}\Bigr)\dd s.
\]
The integrand lies in $L^2(0,T)$, because $Cy\in C([0,T];H)$, $y,p\in C([0,T];H)$ and $y_d\in L^2(0,T;H)$. This gives $\sigma_u\in H^1(0,T)$ and \eqref{eq:switching}. If $y_d$ is continuous, so is the integrand, and $\sigma_u\in C^1([0,T])$.
\end{proof}

The identity separates the two sources of temporal variation of $\sigma$: the commutator, which measures how far the control operator is from commuting with the free dynamics, and the running tracking term.

\begin{remark}\label{rem:lie}
For a control-affine system $x'=f(x)+u\,g(x)$ in $\R^n$ with a terminal cost, the switching function $\ip{g(x)}{p}$ satisfies $\frac{\mathrm d}{\mathrm dt}\ip{g(x)}{p}=\ip{[f,g](x)}{p}$ along extremals, where $[f,g]:=g'f-f'g$ is the Lie bracket \citep{SchattlerLedzewicz2012}. For $f(y)=Ay$ and $g(y)=By$ one has $[f,g](y)=BAy-ABy=Cy$. Identity \eqref{eq:switching} is therefore the infinite-dimensional version of this formula, with an additional term from the running cost. The conditions $y_0\in D(A)$ and $B(D(A))\subset D(A)$ keep the state in $D(A)$, and this is what gives $Cy$ a meaning when $A$ is unbounded; in Example~\ref{ex:heat} below, $C$ is a first-order differential operator. In infinite dimensions, the equation for $By$ derived in Step~2 of the proof also appears in \citet[Section~2]{AronnaBonnansKroner2018}, under an assumption on the extension of the commutator, where it serves the Goh transform for problems that are linear in the control.
\end{remark}

If $A$ and $B$ commute, the commutator vanishes, and the switching function cannot vary when only the final state is penalized.

\begin{corollary}\label{cor:commuting}
Suppose $BS(t)=S(t)B$ for all $t\ge0$ and $\gamma_2=0$. Then $\sigma_u$ is constant on $[0,T]$ for every $u\in\Uad$. Consequently, for every $\beta\ge0$, every control $\bar u\in\Uad$ satisfying \eqref{eq:law} is constant in time, and so is every control satisfying the pointwise law of Remark~\ref{rem:general}. It is either identically zero or nonzero on all of $[0,T]$.
\end{corollary}

\begin{proof}
For $x\in D(A)$, $(S(h)Bx-Bx)/h=B(S(h)x-x)/h\to BAx$, so $B(D(A))\subset D(A)$ and $C=0$. If $y_0\in D(A)$, Theorem~\ref{thm:switching} gives $\sigma_u'=0$. For general $y_0$, take $y_0^k\in D(A)$ converging to $y_0$. For fixed $u$, the state is linear in $y_0$, with $\norm{y_u}_C\le M_TE_T\norm{y_0}$ by Lemma~\ref{lem:wellposed}, and the adjoint is affine in the state; so the state, the adjoint and $\sigma_u$ depend continuously on $y_0$, uniformly on $[0,T]$. So $\sigma_u$ is a uniform limit of constants and is constant. Then $\bar u=\Phi(\bar\sigma)$ is constant by \eqref{eq:law}.
\end{proof}

The case $B=I$ shows this directly.

\begin{example}\label{ex:identity}
If $By=y$, then $y_u(t)=\exp\bigl(\int_0^tu\dd s\bigr)S(t)y_0$, so for $\gamma_2=0$ the tracking term depends on $u$ only through $\int_0^Tu\dd t$. Among controls with a given integral, the constant one minimizes $\frac{\alpha}{2}\norm{u}^2+\beta\norm{u}_{L^1}$, in line with Corollary~\ref{cor:commuting}; the same holds whenever $B$ commutes with $S(t)$, since then $y_u(t)=e^{U(t)B}S(t)y_0$ with $U(t):=\int_0^tu\dd s$. For $\alpha=0$, every control of constant sign with the optimal integral is a solution, so uniqueness fails for every horizon. This is why we keep $\alpha>0$.
\end{example}

When the switching function crosses the thresholds $\pm\beta$ transversally, the support of a control satisfying \eqref{eq:law} consists of finitely many intervals.

\begin{corollary}\label{cor:finite}
Let $y_d\in C([0,T];H)$, $\beta>0$ and $a<0<b$, and let $\bar u\in\Uad$ satisfy \eqref{eq:law}. If $\bar\sigma'(t)\neq0$ at every $t\in[0,T]$ where $\abs{\bar\sigma(t)}=\beta$, then $\supp\bar u$ is a finite union of intervals. Moreover, $\bar u$ is Lipschitz continuous on $[0,T]$; this part does not need the transversality assumption.
\end{corollary}

\begin{proof}
By Theorem~\ref{thm:switching}, $\bar\sigma\in C^1([0,T])$. Each zero of $\abs{\bar\sigma}-\beta$ is a point where $\bar\sigma=\pm\beta$ and $\bar\sigma'\neq0$, hence it is isolated, and $[0,T]$ contains finitely many. Between consecutive zeros $\abs{\bar\sigma}-\beta$ has constant sign, so $\{\bar u\neq0\}=\{\abs{\bar\sigma}>\beta\}$ is a finite union of intervals. Finally $\bar u=\Phi(\bar\sigma)$ is Lipschitz, as $\Phi$ is Lipschitz and $\bar\sigma$ is $C^1$.
\end{proof}

For the bilinear heat equation used in Section~\ref{sec:numerics}, multiplication by a smooth coefficient leaves $D(A)$ invariant, and the commutator can be computed explicitly.

\begin{example}\label{ex:heat}
Let $\Omega\subset\R^n$ be a bounded domain with $C^2$ boundary, $H=L^2(\Omega)$, $A=\nu\Delta$ with $D(A)=H^2(\Omega)\cap H_0^1(\Omega)$, and $(By)(x):=m(x)y(x)$ with $m\in W^{2,\infty}(\Omega)$. Then $B(D(A))\subset D(A)$, and
\[
Cy=m\,\nu\Delta y-\nu\Delta(my)=-\nu\bigl(2\nabla m\cdot\nabla y+(\Delta m)\,y\bigr).
\]
Identity \eqref{eq:switching} becomes
\[
\sigma'=-\nu\int_\Omega\bigl(2\nabla m\cdot\nabla y+(\Delta m)\,y\bigr)p\dd x-\gamma_2\int_\Omega m\,y\,(y-y_d)\dd x.
\]
A constant $m$ gives $C=0$ and falls under Corollary~\ref{cor:commuting}. A spatially varying $m$ makes $\sigma$ vary in time through $\nabla m$ and $\Delta m$; Section~\ref{sec:numerics} uses such an $m$.
\end{example}

\section{Uniqueness and the role of the supports}\label{sec:uniqueness}

We now turn to uniqueness. Two distinct controls in $\Uad$ satisfying \eqref{eq:law} must differ on a set of measure at least $\alpha/L_\sigma$, and this set lies in the union of their supports. With the $L^1$ term the supports are the sets where the switching function exceeds $\beta$, so short supports force uniqueness; for exponentially stable systems steered to zero, this yields a unique optimal control for every horizon.

We first extend Lemma~\ref{lem:lipschitz} to the adjoint state and the switching function. Estimates of this type underlie the classical uniqueness proofs for bilinear problems \citep{BradleyLenhart1997,BouhamedElKaboussBouzahir2022}; here the control is measured in the $L^1$ norm.

\begin{lemma}\label{lem:sigmalip}
For $u_1,u_2\in\Uad$, $\norm{p_{u_1}-p_{u_2}}_C\le L_p\norm{u_1-u_2}_{L^1}$ and $\norm{\sigma_{u_1}-\sigma_{u_2}}_C\le L_\sigma\norm{u_1-u_2}_{L^1}$, where
\begin{equation}\label{eq:Lsigma}
L_p:=M_TE_T\bigl(\gamma_1L_y+\gamma_2TL_y+\norm{B}K_p\bigr),\qquad L_\sigma:=\norm{B}\bigl(L_yK_p+K_yL_p\bigr).
\end{equation}
\end{lemma}

\begin{proof}
The difference $\delta p:=p_{u_1}-p_{u_2}$ satisfies
\[
\delta p(t)=S(T-t)^*\gamma_1\delta y(T)+\int_t^TS(s-t)^*\bigl[u_1B^*\delta p+(u_1-u_2)B^*p_{u_2}+\gamma_2\delta y\bigr]\dd s,
\]
with $\delta y:=y_{u_1}-y_{u_2}$. By Lemma~\ref{lem:lipschitz}, $\norm{\delta p(t)}\le M_T(\gamma_1L_y+\norm{B}K_p+\gamma_2TL_y)\norm{u_1-u_2}_{L^1}+M_T\norm{B}\int_t^T\abs{u_1}\norm{\delta p}\dd s$, and Gronwall's lemma gives the first estimate. The second follows from $\sigma_{u_1}-\sigma_{u_2}=\ip{B\delta y}{p_{u_1}}+\ip{By_{u_2}}{\delta p}$.
\end{proof}

We now prove the separation estimate announced at the beginning of this section.

\begin{theorem}\label{thm:separation}
Let $\bar u_1\neq\bar u_2$ be controls in $\Uad$ satisfying \eqref{eq:law}, with switching functions $\bar\sigma_1,\bar\sigma_2$, and let $D:=\{t\in[0,T]\colon\bar u_1(t)\neq\bar u_2(t)\}$. Then
\begin{equation}\label{eq:separation}
\abs{D}\ge\frac{\alpha}{L_\sigma}.
\end{equation}
Since $D\subset\supp\bar u_1\cup\supp\bar u_2$, in particular $\abs{\supp\bar u_1}+\abs{\supp\bar u_2}\ge\alpha/L_\sigma$. If $\beta>0$ and $a<0<b$, then $\supp\bar u_i=\{\abs{\bar\sigma_i}>\beta\}$ for $i=1,2$.
\end{theorem}

\begin{proof}
By Corollary~\ref{cor:structure}(2), $\abs{\bar u_1(t)-\bar u_2(t)}\le\abs{\bar\sigma_1(t)-\bar\sigma_2(t)}/\alpha$ for all $t$, and $\bar u_1-\bar u_2=0$ outside $D$. With Lemma~\ref{lem:sigmalip},
\[
\norm{\bar u_1-\bar u_2}_{L^1}=\int_D\abs{\bar u_1-\bar u_2}\dd t\le\frac{\abs{D}}{\alpha}\norm{\bar\sigma_1-\bar\sigma_2}_C\le\frac{L_\sigma\abs{D}}{\alpha}\norm{\bar u_1-\bar u_2}_{L^1}.
\]
The controls are continuous and distinct, so $\norm{\bar u_1-\bar u_2}_{L^1}>0$, and \eqref{eq:separation} follows. Both controls vanish outside $\supp\bar u_1\cup\supp\bar u_2$, so $D$ lies in this union, and Corollary~\ref{cor:structure}(4) identifies the supports.
\end{proof}

The classical uniqueness argument integrates the pointwise estimate $\abs{\bar u_1-\bar u_2}\le\abs{\bar\sigma_1-\bar\sigma_2}/\alpha$ over the whole interval $[0,T]$, which yields a contraction only for short horizons (Remark~\ref{rem:classical}). The new point in Theorem~\ref{thm:separation} is to integrate only over the set $D$ where the two controls differ, and to measure their difference in the $L^1$ norm, in which Lemma~\ref{lem:sigmalip} closes the estimate; the length of the horizon is thereby replaced by $\abs{D}$. The proof uses only that the pointwise law is $(1/\alpha)$-Lipschitz, so \eqref{eq:separation} holds for smooth costs as well. What the $L^1$ term changes is the size of the supports. For $\beta=0$ they are the sets $\{\bar\sigma_i\neq0\}$, generically all of $[0,T]$, and \eqref{eq:separation} then only yields the classical small-horizon uniqueness (Remark~\ref{rem:classical}). For $\beta>0$ they are the sets where the switching function exceeds the price $\beta$, and they can be short even when the horizon is long.

\begin{remark}\label{rem:classical}
Since $\abs{D}\le T$, \eqref{eq:separation} cannot hold when $T<\alpha/L_\sigma$. The constant $L_\sigma$ stays bounded as $T\downarrow0$, so this recovers the known uniqueness for short horizons \citep{BradleyLenhart1997,BouhamedElKaboussBouzahir2022}, for $\beta=0$ and $\beta>0$ alike. We do not pursue this case further.
\end{remark}

For $\beta>0$, the separation estimate turns bounds on the supports into uniqueness.

\begin{corollary}\label{cor:shortsupport}
Let $\beta>0$ and $a<0<b$.
\begin{enumerate}[label=\textup{(\arabic*)},leftmargin=*]
\item If every $\bar u\in\Uad$ satisfying \eqref{eq:law} has $\abs{\supp\bar u}<\alpha/(2L_\sigma)$, then exactly one control in $\Uad$ satisfies \eqref{eq:law}, and it is the unique solution of $(\mathrm{P}_\beta)$.
\item If $\bar u\in\Uad$ satisfies \eqref{eq:law} and $\abs{\supp\bar u}=s$, every other control $\tilde u\in\Uad$ satisfying \eqref{eq:law} has $\abs{\supp\tilde u}\ge\alpha/L_\sigma-s$.
\item If $\beta\ge\bar\beta$, the unique solution is $\bar u\equiv0$ \textup{(}Proposition~\textup{\ref{prop:critical})}.
\end{enumerate}
\end{corollary}

\begin{proof}
For two distinct controls satisfying \eqref{eq:law}, Theorem~\ref{thm:separation} gives $\abs{\supp\bar u_1}+\abs{\supp\bar u_2}\ge\abs{D}\ge\alpha/L_\sigma$. This contradicts the assumption of item (1) and proves item (2). A solution exists by Theorem~\ref{thm:existence} and satisfies \eqref{eq:law} by Theorem~\ref{thm:os}.
\end{proof}

Item (2) is an a posteriori test: a computed control satisfying \eqref{eq:law} with short support can only be accompanied by other solutions of \eqref{eq:law} with long supports. To use item (1) one needs an a priori bound on supports. Since $\supp\bar u\subset\{\norm{B}\norm{\bar y(t)}\norm{\bar p(t)}>\beta\}$, any decay of state or adjoint gives such a bound. The next theorem obtains it with constants independent of $T$.

\begin{theorem}\label{thm:longhorizon}
Assume $\norm{S(t)}\le Me^{-\lambda t}$ for $t\ge0$, with $M\ge1$ and $\lambda>0$, and let $y_T=0$, $y_d=0$ and $\beta>0$. Define, independently of $T$ and of the control bounds $a$, $b$,
\begin{equation}\label{eq:stableconstants}
\begin{aligned}
J_0&:=\frac{M^2\norm{y_0}^2}{2}\Bigl(\gamma_1+\frac{\gamma_2}{2\lambda}\Bigr),\quad R:=\frac{J_0}{\beta},\quad \theta:=e^{M\norm{B}R},\\
K&:=M\norm{y_0}\theta,\quad K':=MK\Bigl(\gamma_1+\frac{\gamma_2}{2\lambda}\Bigr)\theta,\quad t_\beta:=\frac{1}{2\lambda}\ln^+\frac{\norm{B}KK'}{\beta},\\
L_1&:=M\norm{B}K\theta,\quad L_2:=M\Bigl(\gamma_1L_1+\norm{B}K'+\frac{\gamma_2L_1}{2\lambda}\Bigr)\theta,\quad L_\infty:=\norm{B}\bigl(L_1K'+KL_2\bigr).
\end{aligned}
\end{equation}
Then, for every $T>0$:
\begin{enumerate}[label=\textup{(\arabic*)},leftmargin=*]
\item every solution $\bar u$ of $(\mathrm{P}_\beta)$ satisfies $\norm{\bar u}_{L^1}\le R$ and $\bar u(t)=0$ for $t\in[t_\beta,T]$;
\item if $\alpha>L_\infty\bigl(1-e^{-2\lambda t_\beta}\bigr)/(2\lambda)$, then $(\mathrm{P}_\beta)$ has a unique solution.
\end{enumerate}
\end{theorem}

\begin{proof}
The proof refines the argument of Theorem~\ref{thm:separation}. The $L^1$ term gives a bound on $\norm{\bar u}_{L^1}$ that does not depend on $T$ (Step~1). With this bound, the estimates of Lemma~\ref{lem:sigmalip} can be repeated with the weight $e^{\lambda t}$, which makes all constants independent of $T$; since $\sigma=\ip{By}{p}$ is the product of two decaying factors, the switching function decays like $e^{-2\lambda t}$ (Steps~2 and~4). This decay confines the supports to $[0,t_\beta]$ (Step~3), and integrating the weighted estimate over this interval closes the argument (Step~5).

\emph{Step 1: an $L^1$ bound independent of $T$.} Since $0\in\Uad$ and all terms of $J$ are nonnegative,
\[
\beta\norm{\bar u}_{L^1}\le J(\bar u)\le J(0)=\frac{\gamma_1}{2}\norm{S(T)y_0}^2+\frac{\gamma_2}{2}\int_0^T\norm{S(t)y_0}^2\dd t\le J_0,
\]
because $\norm{S(t)y_0}\le Me^{-\lambda t}\norm{y_0}$ and $\int_0^Te^{-2\lambda t}\dd t\le1/(2\lambda)$. So $\norm{\bar u}_{L^1}\le R$.

\emph{Step 2: weighted bounds.} Let $u\in\Uad$ with $\norm{u}_{L^1}\le R$, and write $y:=y_u$ and $p:=p_u$. Multiplying \eqref{eq:mild} by $e^{\lambda t}$ and using $\norm{S(t-s)}\le Me^{-\lambda(t-s)}$ gives
\[
e^{\lambda t}\norm{y(t)}\le M\norm{y_0}+M\norm{B}\int_0^t\abs{u(s)}\,e^{\lambda s}\norm{y(s)}\dd s,
\]
and Gronwall's lemma yields $e^{\lambda t}\norm{y(t)}\le M\norm{y_0}e^{M\norm{B}\norm{u}_{L^1}}\le K$. Since $y_T=0$ and $y_d=0$, the adjoint is
\[
p(t)=\gamma_1S(T-t)^*y(T)+\int_t^TS(s-t)^*\bigl[u(s)B^*p(s)+\gamma_2y(s)\bigr]\dd s.
\]
For $s\ge t$ we have $e^{\lambda t}e^{-\lambda(s-t)}e^{-\lambda s}=e^{-2\lambda(s-t)}\le1$. Hence, by the bound on $y$, the function $\chi(t):=e^{\lambda t}\norm{p(t)}$ satisfies
\begin{align*}
\chi(t)&\le\gamma_1MK+\gamma_2MK\int_t^Te^{-2\lambda(s-t)}\dd s+M\norm{B}\int_t^T\abs{u(s)}\chi(s)\dd s\\
&\le MK\Bigl(\gamma_1+\frac{\gamma_2}{2\lambda}\Bigr)+M\norm{B}\int_t^T\abs{u(s)}\chi(s)\dd s,
\end{align*}
and the backward Gronwall lemma gives $\chi\le K'$. Therefore $\norm{y(t)}\le Ke^{-\lambda t}$, $\norm{p(t)}\le K'e^{-\lambda t}$ and
\begin{equation}\label{eq:sigmadecay}
\abs{\sigma_u(t)}\le\norm{B}KK'e^{-2\lambda t},\qquad t\in[0,T].
\end{equation}

\emph{Step 3: localization.} By Steps~1 and~2, every solution satisfies $\abs{\bar\sigma(t)}\le\norm{B}KK'e^{-2\lambda t}\le\beta$ for $t\ge t_\beta$, by the definition of $t_\beta$, so $\bar u(t)=0$ there by Corollary~\ref{cor:structure}(1). This proves item~(1).

\emph{Step 4: weighted Lipschitz bounds.} Let $\bar u_1,\bar u_2$ be solutions, with states $y_i$, adjoints $p_i$ and switching functions $\sigma_i$, and set $\delta u:=\bar u_1-\bar u_2$, $\delta y:=y_1-y_2$, $\delta p:=p_1-p_2$ and $\delta\sigma:=\sigma_1-\sigma_2$. Since $\delta y(t)=\int_0^tS(t-s)[\bar u_1B\delta y+\delta u\,By_2]\dd s$ and $e^{\lambda s}\norm{y_2(s)}\le K$, the argument of Step~2 gives
\[
e^{\lambda t}\norm{\delta y(t)}\le M\norm{B}K\norm{\delta u}_{L^1}+M\norm{B}\int_0^t\abs{\bar u_1(s)}\,e^{\lambda s}\norm{\delta y(s)}\dd s,
\]
hence $e^{\lambda t}\norm{\delta y(t)}\le L_1\norm{\delta u}_{L^1}$ by Gronwall's lemma. The difference of the adjoints satisfies
\[
\delta p(t)=\gamma_1S(T-t)^*\delta y(T)+\int_t^TS(s-t)^*\bigl[\bar u_1B^*\delta p+\delta u\,B^*p_2+\gamma_2\delta y\bigr]\dd s.
\]
With $e^{\lambda s}\norm{\delta y(s)}\le L_1\norm{\delta u}_{L^1}$ and $e^{\lambda s}\norm{p_2(s)}\le K'$, the three source terms contribute at most $\gamma_1ML_1$, $M\norm{B}K'$ and $\gamma_2ML_1/(2\lambda)$ to $e^{\lambda t}\norm{\delta p(t)}$, each times $\norm{\delta u}_{L^1}$, and the backward Gronwall lemma gives $e^{\lambda t}\norm{\delta p(t)}\le L_2\norm{\delta u}_{L^1}$. From $\delta\sigma=\ip{B\delta y}{p_1}+\ip{By_2}{\delta p}$ we obtain
\begin{equation}\label{eq:deltasigma}
\abs{\delta\sigma(t)}\le\norm{B}\bigl(L_1K'+KL_2\bigr)e^{-2\lambda t}\norm{\delta u}_{L^1}=L_\infty e^{-2\lambda t}\norm{\delta u}_{L^1},\qquad t\in[0,T].
\end{equation}

\emph{Step 5: conclusion.} By item~(1), $\delta u=0$ on $[t_\beta,T]$. Both controls satisfy \eqref{eq:law} by Theorem~\ref{thm:os}, so Corollary~\ref{cor:structure}(2) and \eqref{eq:deltasigma} give $\abs{\delta u(t)}\le\abs{\delta\sigma(t)}/\alpha\le(L_\infty/\alpha)e^{-2\lambda t}\norm{\delta u}_{L^1}$. Integrating over $[0,\min\{t_\beta,T\}]$ gives
\[
\norm{\delta u}_{L^1}\le\frac{L_\infty}{\alpha}\int_0^{t_\beta}e^{-2\lambda t}\dd t\;\norm{\delta u}_{L^1}=\frac{L_\infty\bigl(1-e^{-2\lambda t_\beta}\bigr)}{2\lambda\alpha}\norm{\delta u}_{L^1}.
\]
Under the condition of item~(2) the factor is smaller than $1$, so $\delta u=0$.
\end{proof}

The proof uses the optimality of $\bar u_1$ and $\bar u_2$ only through $J(\bar u_i)\le J(0)$ in Step~1, and Steps~3 and~5 use only \eqref{eq:law}. Under the condition of item~(2), every $\tilde u\in\Uad$ that satisfies \eqref{eq:law} and $J(\tilde u)\le J(0)$ therefore coincides with the solution of $(\mathrm{P}_\beta)$; this applies, for instance, to the limit of a descent method started at $u=0$ when this limit satisfies \eqref{eq:law}.

\begin{remark}\label{rem:whyL1}
Two ingredients of Theorem~\ref{thm:longhorizon} come from $\beta>0$: the bound $\norm{\bar u}_{L^1}\le J_0/\beta$, uniform in $T$, and the localization of supports in $[0,t_\beta]$. For $\beta=0$ the cost only gives $\norm{\bar u}_{L^2}\le(2J_0/\alpha)^{1/2}$, hence $\norm{\bar u}_{L^1}\le(2TJ_0/\alpha)^{1/2}$. This grows with $T$ and enters the constants through $\exp(M\norm{B}\norm{\bar u}_{L^1})$, so the argument degenerates as $T\to\infty$.

The control bounds can take the place of the $L^1$ term only when they are small compared with the decay rate. If $\mu:=M\norm{B}\rho<\lambda$, every admissible control preserves exponential stability, $\norm{y_u(t)}\le M\norm{y_0}e^{-(\lambda-\mu)t}$. Steps~2, 4 and~5 of the proof, with Gronwall's lemma replaced by weighted supremum estimates with a rate $\omega\in(0,\lambda-\mu)$ and with the integral in Step~5 taken over $[0,T]$, then show that exactly one control in $\Uad$ satisfies \eqref{eq:law}, for every $T>0$ and every $\beta\ge0$, provided $\alpha$ exceeds a constant independent of $T$. Theorem~\ref{thm:longhorizon} needs no relation between $\rho$ and $\lambda$, and by Remark~\ref{rem:unconstrained} it holds without control bounds. In the experiments of Section~\ref{sec:numerics}, $\mu=6$ and $\lambda\approx0.49$.
\end{remark}

\begin{remark}\label{rem:nonvoid}
The right-hand side of the condition in item~(2) does not depend on $\alpha$. It decreases as $\beta$ grows and vanishes once $\beta\ge\norm{B}KK'$, so for every $\alpha>0$ there is $\beta^*(\alpha)$, independent of $T$, such that the optimal control is unique for all $\beta\ge\beta^*(\alpha)$ and all $T>0$. Conversely, for fixed $\beta>0$ the condition holds for all $\alpha$ above a threshold $\alpha^*(\beta)$ independent of $T$. If moreover $\beta<\beta_0$, the critical weight of Proposition~\ref{prop:critical} for the horizon $T$, the zero control does not satisfy \eqref{eq:law}, whatever $\alpha$. In this regime the theorem certifies a unique optimal control that is not zero.
\end{remark}

When the target is not zero, uniqueness can fail even for $\beta>0$, as the following scalar example shows.

\begin{example}\label{ex:nonunique}
Let $H=\R$, $A=-1$, $B=1$, $y_0=1$, $\gamma_1=1$, $\gamma_2=0$, $y_T=1$, $T=5$, $\alpha=0.1$ and $[a,b]=[-2,2]$; the target is not zero, so Theorem~\ref{thm:longhorizon} does not apply. Since $A$ and $B$ commute, every control satisfying \eqref{eq:law} is a constant $u=c/T$ (Corollary~\ref{cor:commuting}), and $J$ reduces to
\[
f(c)=\tfrac12\bigl(e^{c-5}-1\bigr)^2+\frac{\alpha c^2}{2T}+\beta\abs{c}.
\]
Its critical points are listed in Table~\ref{tab:nonunique}. The zero control satisfies \eqref{eq:law} for $\beta\ge\beta_0=e^{-5}(1-e^{-5})\approx0.0067$, two global solutions coexist at $\beta\approx0.0513$, and uniqueness returns only for $\beta>0.164$. The two local minimizers differ on all of $[0,T]$, which \eqref{eq:separation} allows: the $L^1$ term does not convexify the problem, and the uniqueness in Theorem~\ref{thm:longhorizon} comes from the decay of state and adjoint, not from convexity.
\end{example}

\begin{table}[t]
\centering
\caption{Controls satisfying \eqref{eq:law} in Example~\ref{ex:nonunique} (value of the constant control).}\label{tab:nonunique}
\small
\begin{tabular}{@{}lll@{}}
\toprule
$\beta$ & Solutions of \eqref{eq:law} & Global solution\\
\midrule
$0$ & $0.123$ (local min), $0.315$ (local max), $0.977$ (local min) & $0.977$\\
$0.01$ & $0$ (local min), $0.420$ (local max), $0.974$ (local min) & $0.974$\\
$0.05$ & $0$, $0.581$, $0.961$ & $0.961$\\
$0.06$ & $0$, $0.607$, $0.957$ & $0$\\
$0.10$ & $0$, $0.694$, $0.939$ & $0$\\
$0.20$ & $0$ & $0$\\
\bottomrule
\end{tabular}
\end{table}

\section{Second-order sufficient conditions}\label{sec:ssc}

A control satisfying \eqref{eq:law} at which $F''$ is positive on the critical cone is a strict local minimizer with quadratic growth. The $L^1$ term shrinks the critical cone to the support, so the sufficient condition of Corollary~\ref{cor:ssc}, which can be evaluated on a computed solution, asks $\alpha$ to dominate the measure of the support rather than the length of the horizon. Throughout this section $a<0<b$.

We first compute the second derivative of the smooth part of the cost.

\begin{proposition}\label{prop:second}
For $u\in\Uad$ and $v\in L^2(0,T)$, with $z_v$ from \eqref{eq:linearized},
\begin{equation}\label{eq:second}
F''(u)[v,v]=\alpha\norm{v}_{L^2}^2+\gamma_1\norm{z_v(T)}^2+\gamma_2\norm{z_v}_{L^2(0,T;H)}^2+2\int_0^Tv(t)\,\ip{Bz_v(t)}{p_u(t)}\dd t.
\end{equation}
\end{proposition}

\begin{proof}
The same expansion, for costs that are linear in the control, is derived in \citet[Section~3]{AronnaBonnansKroner2018}. Write $F(u)=\Psi(G(u))+\frac{\alpha}{2}\norm{u}^2$, where $\Psi$ is the tracking part. Then $F''(u)[v,v]=\Psi''(y_u)[z_v,z_v]+\Psi'(y_u)w_v+\alpha\norm{v}^2$, and the first term equals $\gamma_1\norm{z_v(T)}^2+\gamma_2\norm{z_v}^2$. For the second, apply Lemma~\ref{lem:duality} with $t_0=0$, $z=w_v$, $f=uBw_v+2vBz_v$ and $q=p_u$. As in Proposition~\ref{prop:derivative}, the terms containing $u$ cancel and $\Psi'(y_u)w_v=2\int_0^Tv\ip{Bz_v}{p_u}\dd t$.
\end{proof}

To state the sufficient condition, let $\bar u\in\Uad$ satisfy \eqref{eq:law}, let $\bar\sigma$ be its switching function and $\bar d:=\alpha\bar u+\bar\sigma$. Following \citet{CasasHerzogWachsmuth2012}, we define the critical cone $C_{\bar u}$ as the set of $v\in L^2(0,T)$ such that, for a.e.\ $t$:
\begin{itemize}[leftmargin=*]
\item $v(t)=0$ if $\bar u(t)=0$ and $\abs{\bar\sigma(t)}<\beta$, and $v(t)\bar\sigma(t)\le0$ if $\bar u(t)=0$ and $\abs{\bar\sigma(t)}=\beta$;
\item $v(t)=0$ if $\bar u(t)=a$ and $\bar d(t)>\beta$, and $v(t)\ge0$ if $\bar u(t)=a$ and $\bar d(t)=\beta$;
\item $v(t)=0$ if $\bar u(t)=b$ and $\bar d(t)<-\beta$, and $v(t)\le0$ if $\bar u(t)=b$ and $\bar d(t)=-\beta$.
\end{itemize}
Equivalently, $C_{\bar u}$ is the set of directions $h$ tangent to $\Uad$ at $\bar u$ with $\ell(h)=0$, where $\ell(h):=F'(\bar u)h+\beta j'(\bar u;h)$ and $j'(\bar u;h)=\int_{\{\bar u>0\}}h-\int_{\{\bar u<0\}}h+\int_{\{\bar u=0\}}\abs{h}$. Every $v\in C_{\bar u}$ vanishes outside
\begin{equation}\label{eq:Sbar}
\bar S:=\{\bar u\notin\{a,0,b\}\}\cup\{\bar u=0,\ \abs{\bar\sigma}=\beta\}\cup\{\bar u=a,\ \bar d=\beta\}\cup\{\bar u=b,\ \bar d=-\beta\}.
\end{equation}

With this notation, the second-order sufficient condition reads as follows.

\begin{theorem}\label{thm:ssc}
Let $\bar u\in\Uad$ satisfy \eqref{eq:law} and assume
\begin{equation}\label{eq:ssc}
F''(\bar u)[v,v]>0\qquad\text{for all }v\in C_{\bar u}\setminus\{0\}.
\end{equation}
Then there are $\varepsilon>0$ and $\delta>0$ such that $J(u)\ge J(\bar u)+\frac{\delta}{2}\norm{u-\bar u}_{L^2}^2$ for all $u\in\Uad$ with $\norm{u-\bar u}_{L^2}\le\varepsilon$. In particular, $\bar u$ is a strict local solution of $(\mathrm{P}_\beta)$.
\end{theorem}

\begin{proof}
The proof is the classical contradiction argument for sparse control problems \citep[Section~3]{CasasHerzogWachsmuth2012}, see also \citet{CasasTroltzsch2012}, and we only indicate its steps. If the claim fails, there are $u_k\in\Uad$ with $\rho_k:=\norm{u_k-\bar u}_{L^2}\in(0,1/k)$ and $J(u_k)<J(\bar u)+\rho_k^2/(2k)$, and, for a subsequence, $v_k:=(u_k-\bar u)/\rho_k\rightharpoonup v$ in $L^2(0,T)$. The functional $\ell$ is convex and continuous, it is nonnegative on the closed convex cone $\mathcal T:=\{h\colon h\ge0\text{ on }\{\bar u=a\},\ h\le0\text{ on }\{\bar u=b\}\}$ by \eqref{eq:vi}, and $\mathcal T$ contains all $v_k$, hence $v$. Since $j(u_k)-j(\bar u)\ge\rho_kj'(\bar u;v_k)$ by convexity, Taylor expansions of $F$ of first and second order give $\limsup_k\ell(v_k)\le0$ and $\limsup_kF''(\bar u)[v_k,v_k]\le0$. The first inequality and the weak lower semicontinuity of $\ell$ give $\ell(v)=0$, that is, $v\in C_{\bar u}$.

The only point specific to the bilinear problem is the weak sequential continuity of the form $Q(h):=F''(\bar u)[h,h]-\alpha\norm{h}^2$, which collects the last three terms of \eqref{eq:second}. If $h_k\rightharpoonup h$ in $L^2(0,T)$, then $\zeta_k:=z_{h_k}-z_h$ is the mild solution of $\zeta'=A\zeta+\bar uB\zeta+(h_k-h)B\bar y$, $\zeta(0)=0$, and Gronwall's lemma together with Lemma~\ref{lem:bms}, applied with $u=\bar u$ and $w_k=h_k-h$, gives $z_{h_k}\to z_h$ in $C([0,T];H)$. As $\ip{Bz_{h_k}}{\bar p}\to\ip{Bz_h}{\bar p}$ uniformly, $Q(h_k)\to Q(h)$. Hence, since $\norm{v_k}=1\ge\norm{v}$,
\[
0\ge\lim_{k\to\infty}F''(\bar u)[v_k,v_k]=\alpha+Q(v)\ge\alpha\norm{v}^2+Q(v)=F''(\bar u)[v,v].
\]
By \eqref{eq:ssc}, $v=0$, and then $\alpha\le0$, a contradiction.
\end{proof}

The condition can be verified through an explicit inequality that involves only the computed solution and the measure of $\bar S$.

\begin{corollary}\label{cor:ssc}
Let $\bar u\in\Uad$ satisfy \eqref{eq:law}, with state $\bar y$ and adjoint $\bar p$, let $\bar S$ be as in \eqref{eq:Sbar}, and let $\kappa:=M_T\norm{B}$. If
\begin{equation}\label{eq:ssctest}
\alpha>c(\bar u)\,\abs{\bar S},\qquad c(\bar u):=2\kappa\norm{B}\,\norm{\bar y}_C\,\norm{\bar p}_C\,e^{\kappa\norm{\bar u}_{L^1(0,T)}},
\end{equation}
then \eqref{eq:ssc} holds, and $\bar u$ is a strict local solution with quadratic growth. Since $\norm{\bar y}_C\le K_y$, $\norm{\bar p}_C\le K_p$ and $e^{\kappa\norm{\bar u}_{L^1}}\le E_T$, we have $c(\bar u)\le2\norm{B}K_pL_y$; in particular \eqref{eq:ssctest} holds if $\alpha>2\norm{B}K_pL_y\abs{\bar S}$.
\end{corollary}

\begin{proof}
Let $v\in C_{\bar u}\setminus\{0\}$. Then $v=0$ outside $\bar S$, so $\norm{v}_{L^1}\le\abs{\bar S}^{1/2}\norm{v}_{L^2}$. By \eqref{eq:zbound}, $\norm{z_v}_C\le\kappa\norm{\bar y}_Ce^{\kappa\norm{\bar u}_{L^1}}\norm{v}_{L^1}$, so
\[
\Bigl|2\int_0^Tv\ip{Bz_v}{\bar p}\dd t\Bigr|\le2\norm{B}\norm{\bar p}_C\norm{z_v}_C\norm{v}_{L^1}\le c(\bar u)\norm{v}_{L^1}^2\le c(\bar u)\abs{\bar S}\norm{v}_{L^2}^2.
\]
The two tracking terms in \eqref{eq:second} are nonnegative, hence $F''(\bar u)[v,v]\ge(\alpha-c(\bar u)\abs{\bar S})\norm{v}_{L^2}^2>0$.
\end{proof}

\begin{remark}\label{rem:sscL1}
For $\beta=0$, $\bar S$ contains every $t$ with $a<\bar u(t)<b$, which is typically almost all of $[0,T]$, so \eqref{eq:ssctest} asks $\alpha$ to dominate $c(\bar u)\,T$. For $\beta>0$, the strongly inactive set $\{\bar u=0,\ \abs{\bar\sigma}<\beta\}$ is removed from $\bar S$. Under the transversality assumption of Corollary~\ref{cor:finite}, the set $\{\bar u=0,\ \abs{\bar\sigma}=\beta\}$ is finite, and $\abs{\bar S}$ is the measure of the support up to the weakly active sets $\{\bar u=a,\ \bar d=\beta\}$ and $\{\bar u=b,\ \bar d=-\beta\}$. The factor $c(\bar u)$ depends on the horizon only through $M_T$, $\norm{\bar y}_C$, $\norm{\bar p}_C$ and $e^{\kappa\norm{\bar u}_{L^1}}\le e^{\kappa\rho\abs{\supp\bar u}}$; for exponentially stable semigroups $M_T\le M$. The test thus asks $\alpha$ to dominate the time during which the control acts, and it can be evaluated on a computed solution. Section~\ref{sec:numerics} shows that it is far from sharp, while \eqref{eq:ssc} itself can be checked numerically. Quadratic growth is also the usual starting point for Lipschitz stability with respect to the data, which we do not pursue here.
\end{remark}

\section{Numerical experiments}\label{sec:numerics}

We illustrate the results on the bilinear heat equation of Example~\ref{ex:heat} with $\Omega=(0,1)$, $\nu=0.05$, $m(x)=\exp(-((x-0.35)/0.12)^2)$, $y_0(x)=\sin(\pi x)$, $[a,b]=[-6,6]$ and $\alpha=10^{-3}$. Since $y_0\in D(A)$ and $m$ is smooth, the results of Section~\ref{sec:switching} apply. Space is discretized by finite differences with $99$ interior nodes ($h=0.01$) and time by the step $\Delta t=0.0025$, with piecewise constant controls $u_n$ on $(t_{n-1},t_n]$. The state is advanced by the Lie splitting
\begin{equation}\label{eq:scheme}
y_n=E(u_n)\,R\,y_{n-1},\qquad R:=(I-\Delta t\,\nu L_h)^{-1},\qquad E(u):=\diag\bigl(e^{\Delta t\,u\,m(x_i)}\bigr),
\end{equation}
which treats diffusion implicitly, integrates the bilinear reaction exactly and keeps the state positive. When $m$ is constant, $E(u)$ commutes with $R$, so Corollary~\ref{cor:commuting} also holds for the discrete problem. The discrete adjoint and switching function,
\begin{equation}\label{eq:discadjoint}
\begin{gathered}
p_N=\gamma_1(y_N-y_T)+\gamma_2\Delta t\,(y_N-y_{d,N}),\qquad p_n=R\,E(u_{n+1})\,p_{n+1}+\gamma_2\Delta t\,(y_n-y_{d,n}),\\
\sigma_n=h\sum_im(x_i)\,y_{n,i}\,p_{n,i},
\end{gathered}
\end{equation}
give the exact gradient of the discrete tracking functional; finite differences confirm it to a relative error below $2\times10^{-7}$. We minimize with L-BFGS-B \citep{ByrdLuNocedalZhu1995} on the splitting $u=u^+-u^-$, followed by proximal gradient steps $u\leftarrow P\bigl(S_{\tau\beta}(u-\tau\sigma(u))/(1+\tau\alpha)\bigr)$ \citep{BeckTeboulle2009}, and stop when $\max_n\abs{u_n-P(-S_\beta(\sigma_n)/\alpha)}<10^{-9}$. Semismooth Newton methods \citep{Stadler2009,Ulbrich2011} would converge faster but need second derivatives, whereas this method needs only the gradient. Every computed solution was reproduced from several random initial guesses.

\begin{figure}[t]
\centering
\includegraphics[width=0.9\textwidth]{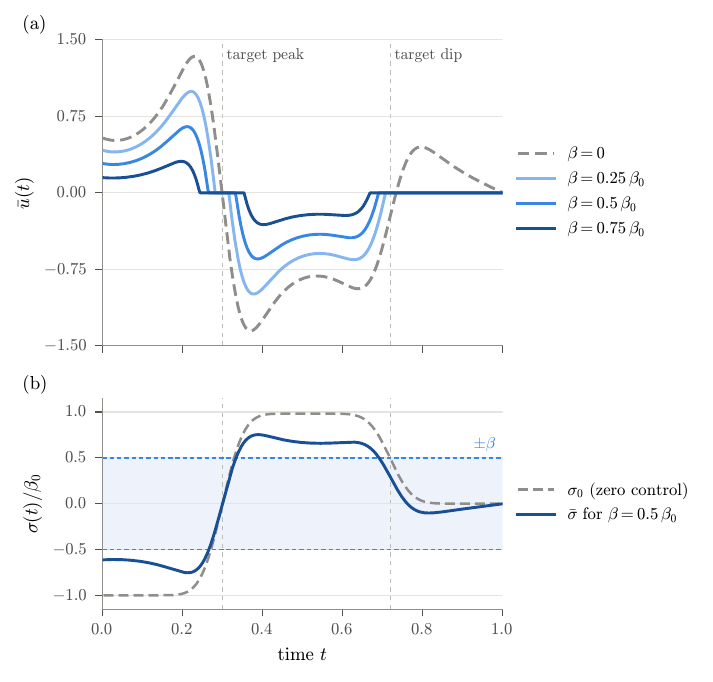}
\caption{Scenario A. (a) Optimal controls for $\beta/\beta_0=0$, $0.25$, $0.5$ and $0.75$. (b) Switching function $\bar\sigma$ for $\beta=0.5\beta_0$ and switching function $\sigma_0$ of the zero control, divided by $\beta_0$; the shaded band is the dead zone $[-\beta,\beta]$, where $\bar u=0$. Dashed vertical lines mark the peak and the dip of the target.}\label{fig:controls}
\end{figure}

\begin{figure}[t]
\centering
\includegraphics[width=0.83\textwidth]{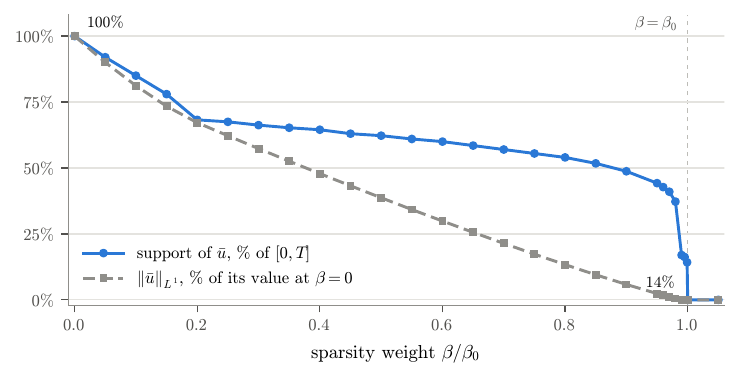}
\caption{Support of the optimal control (percentage of $[0,T]$) and its $L^1$ norm (percentage of its value at $\beta=0$) along the $\beta$-path of scenario A. The control vanishes at the critical weight $\beta_0$ of Proposition~\ref{prop:critical} (dashed line).}\label{fig:betapath}
\end{figure}

In scenario A the target rises by $40\%$ around $t=0.3$ and falls by $30\%$ around $t=0.72$: $y_d(t)=\theta(t)y_{\mathrm{free}}(t)$, where $y_{\mathrm{free}}$ is the uncontrolled state and $\theta(t)=1+0.4e^{-((t-0.3)/0.06)^2}-0.3e^{-((t-0.72)/0.06)^2}$. We take $\gamma_1=0$, $\gamma_2=1$ and $T=1$; the critical weight is $\beta_0=2.565\times10^{-3}$. Without the $L^1$ term the control never vanishes (Figure~\ref{fig:controls}). With it, the control is switched off around its sign change and after $t\approx0.7$, and each off-interval starts where $\bar\sigma$ crosses the dead zone $[-\beta,\beta]$ transversally, as in Corollary~\ref{cor:finite}. Along the $\beta$-path of Figure~\ref{fig:betapath}, the $L^1$ norm of the control decreases to zero, and the control is exactly zero at $\beta=\beta_0$, as Proposition~\ref{prop:critical} predicts. The support, however, still covers $14\%$ of the horizon at $\beta=0.999\beta_0$: on $[0,0.15]$ the target coincides with the free state, the tracking term in \eqref{eq:switching} almost vanishes, and $\sigma_0$ stays within $0.05\%$ of $-\beta_0$. For $\beta=0$ and $\beta=0.5\beta_0$, eight random initial guesses gave the same control to within $1.2\times10^{-9}$. Table~\ref{tab:mesh} repeats the case $\beta=0.5\beta_0$ on five grids: the critical weight changes by less than $0.05\%$, the support intervals agree to within one time step, and the control converges with order about one, as expected for the splitting \eqref{eq:scheme}.

\begin{table}[t]
\centering
\caption{Scenario A with $\beta=1.2825\times10^{-3}$ on five grids. $\beta_0$ is the discrete critical weight, $\bar u_{\mathrm{ref}}$ the control on the finest grid, and the order is $\log_2$ of the ratio of consecutive errors.}\label{tab:mesh}
\footnotesize
\setlength{\tabcolsep}{4pt}
\begin{tabular}{@{}cccccccc@{}}
\toprule
$h$ & $\Delta t$ & $\beta_0\cdot10^{3}$ & $\supp\bar u$ & $\norm{\bar u}_{L^1}$ & $J\cdot10^{3}$ & $\norm{\bar u-\bar u_{\mathrm{ref}}}_{L^1}$ & order\\
\midrule
$0.04$ & $0.01$ & $2.5641$ & $[0,0.260]\cup[0.330,0.690]$ & $0.2632$ & $2.9184$ & $7.5\times10^{-3}$ & --\\
$0.02$ & $0.005$ & $2.5648$ & $[0,0.265]\cup[0.330,0.690]$ & $0.2650$ & $2.9148$ & $3.6\times10^{-3}$ & $1.07$\\
$0.01$ & $0.0025$ & $2.5650$ & $[0,0.265]\cup[0.3325,0.690]$ & $0.2658$ & $2.9134$ & $1.7\times10^{-3}$ & $1.07$\\
$0.005$ & $0.00125$ & $2.5651$ & $[0,0.265]\cup[0.3325,0.690]$ & $0.2662$ & $2.9127$ & $7.9\times10^{-4}$ & $1.11$\\
$0.0025$ & $0.000625$ & $2.5652$ & $[0,0.2644]\cup[0.3325,0.690]$ & $0.2664$ & $2.9124$ & -- & --\\
\bottomrule
\end{tabular}
\end{table}

In scenario B only the final state is penalized ($\gamma_1=1$, $\gamma_2=0$, $y_T=1.3\,y_{\mathrm{free}}(T)$), and $\beta$ is half the critical weight. For $m\equiv1$, where $B=I$ commutes with $A$, the computed control is constant to $1.8\times10^{-9}$ and agrees to $10^{-12}$ with the minimizer of the scalar problem of Example~\ref{ex:identity}, so no temporal sparsity occurs. With the bump $m$, the same cost gives a control that is active on $59\%$ of the horizon.

In scenario C the state is steered to zero ($y_d=0$, $\gamma_1=0$, $\gamma_2=1$, $\beta=0.042$). The heat semigroup decays with rate $\lambda=\nu\pi^2\approx0.49$, so Theorem~\ref{thm:longhorizon} applies. Table~\ref{tab:horizons} shows that the support ends before $t=0.43$ for every horizon, that the controls on $[0,2]$ approach the one for $T=16$ roughly at the rate $e^{-2\lambda T}$ of \eqref{eq:sigmadecay}, and that multi-start runs agree for every horizon. The constants of the theorem are far more pessimistic, $t_\beta\approx21$ and $L_\infty\approx4.6\times10^{10}$, so the theorem itself does not certify uniqueness here.

\begin{table}[t]
\centering
\caption{Scenario C (zero targets): optimal controls for four horizons with the same time step.}\label{tab:horizons}
\small
\begin{tabular}{@{}cccccc@{}}
\toprule
$T$ & End of support & $\norm{\bar u}_{L^1}$ & Optimal cost $J$ & $\max_{[0,2]}\abs{\bar u_T-\bar u_{16}}$ & Multi-start spread\\
\midrule
$2$ & $0.315$ & $1.065$ & $0.1795$ & $1.20$ & $6\times10^{-10}$\\
$4$ & $0.405$ & $1.262$ & $0.1956$ & $0.215$ & $8\times10^{-10}$\\
$8$ & $0.420$ & $1.293$ & $0.1980$ & $4.5\times10^{-3}$ & $1\times10^{-9}$\\
$16$ & $0.4225$ & $1.294$ & $0.1981$ & $0$ & $1\times10^{-10}$\\
\bottomrule
\end{tabular}
\end{table}

Finally, Table~\ref{tab:ssc} checks the second-order condition \eqref{eq:ssc} at the computed solutions. Strict complementarity holds on the grid, so the discrete critical cone consists of the controls supported on the free set $\bar S_h:=\{n\colon a<u_n<b,\ u_n\neq0\}$ (for $\beta=0$ the condition $u_n\neq0$ is dropped). We compute the smallest eigenvalue $\lambda_{\min}$ of the reduced Hessian on this set from central differences of the exact gradient. It exceeds $\alpha$ in every case, so each computed solution is a strict local minimizer of the discrete problem. The a posteriori test \eqref{eq:ssctest}, evaluated with $M_T=1$ and $\norm{B}=1$, fails by a factor between $12$ and $707$; as predicted, its value stays bounded as $T$ grows.

\begin{table}[t]
\centering
\caption{Second-order check at the computed solutions: measure of the free set, smallest eigenvalue of the reduced Hessian on the discrete critical cone, and ratio of the right-hand side of the test \eqref{eq:ssctest} to $\alpha$.}\label{tab:ssc}
\small
\begin{tabular}{@{}lccc@{}}
\toprule
Solution & $\abs{\bar S}$ & $\lambda_{\min}/\alpha$ & $c(\bar u)\abs{\bar S}/\alpha$\\
\midrule
A, $\beta=0$ & $1.0000$ & $1.00006$ & $29.9$\\
A, $\beta=0.25\beta_0$ & $0.6750$ & $1.00010$ & $16.8$\\
A, $\beta=0.5\beta_0$ & $0.6225$ & $1.00011$ & $14.6$\\
A, $\beta=0.75\beta_0$ & $0.5550$ & $1.00011$ & $12.3$\\
C, $T=2$ & $0.2225$ & $1.00042$ & $360$\\
C, $T=4$ & $0.3000$ & $1.00041$ & $640$\\
C, $T=8$ & $0.3150$ & $1.00040$ & $700$\\
C, $T=16$ & $0.3175$ & $1.00040$ & $707$\\
\bottomrule
\end{tabular}
\end{table}

\section{Conclusions}\label{sec:conclusion}

Adding an $L^1$ term to a quadratic control cost leaves the classical theory of bilinear optimal control intact and adds structure that a smooth cost cannot produce: zero intervals, explicit thresholds, supports that stay short on long horizons, and uniqueness for every horizon in stable problems steered to zero, with constants that depend neither on the horizon nor on the control bounds. To our knowledge, this is the first uniqueness result for bilinear optimal control problems that is not local in time. The switching identity \eqref{eq:switching} tells when sparsity in time is possible at all: if $A$ and $B$ commute and only the final state is penalized, no strongly convex control cost produces it.

Several questions remain open. Controls that depend on space and time, or several inputs acting through operators $B_1,\dots,B_k$, require parabolic $L^p$ estimates and vector-valued switching functions, and the commutators $[B_i,A]$ and $[B_i,B_j]$ should then govern the switching structure. The limit $\alpha\to0$ leads to bang-off-bang controls and singular arcs, for which the Goh transform of \citet{AronnaBonnansKroner2018} and second-order analysis of bang-bang bilinear problems \citep{CasasWachsmuthWachsmuth2018} are natural tools. Since the control vanishes after a time independent of the horizon, a direct passage to infinite horizons should be possible, as for sparse stabilization of semilinear equations \citep{CasasKunisch2017} and for bilinear problems with quadratic costs \citep{BreitenKunischPfeiffer2018}. Finally, the constants of Theorem~\ref{thm:longhorizon} and of the test \eqref{eq:ssctest} are far from what the experiments show; sharper estimates for specific dissipative operators, and an extension to unbounded control operators such as boundary bilinear controls, are natural next steps.

\section*{Disclosure statement}
No potential conflict of interest was reported by the author(s).

\section*{Funding}
The authors received no funding for this work.

\section*{Data availability statement}
Data sharing is not applicable to this article as no new data were created or analysed in this study.

\end{document}